\input amstex
\documentstyle{amsppt}
\NoBlackBoxes
\magnification 1200
\hsize 6.5truein
\vsize 8truein
\topmatter
\title Garaev's inequality in Finite Fields not of prime order \endtitle
\author  Nets Hawk Katz and Chun-Yen Shen \endauthor
\affil Indiana University \endaffil
 \thanks 
  The first  author was supported by NSF grant DMS 0432237.
\endthanks 
\endtopmatter

\head \S 0 Introduction \endhead

Let $A$ be a subset of $F=F_{p^k}$, the field of $p^k$ elements with $p$ prime.

We let
$$A+A=\{a+b: a \in A, b \in A \},$$
and
$$AA=\{ab: a \in A, b \in A \}.$$

It is fun (and useful) to prove lower bounds on $\max(|A+A|,|AA|)$ (see e.g. [BKT],[BGK],[G]).
Recently, Garaev [G] showed that when $k=1$ and $|A| < p^{{1 \over 2}}$ one has the estimate
$$\max(|A+A|,|AA|) \gtrapprox |A|^{{15 \over 14}}.$$
The authors in [KS] slightly improved this to
$$\max(|A+A|,|AA|) \gtrapprox |A|^{{14 \over 13}}.$$

In the present paper, we extend Garaev's techniques to the set of fields which are not necessarily of prime order.
Precisely, we prove

\proclaim{Main Theorem} Let $F$ be a finite field. Suppose that $A$ is a subset of $F$
  so that for any $A^{\prime} \subset A$ 
with $|A^{\prime}| \geq |A|^{{47 \over 48}}$ and for
any $G \subset F$
 a subfield (not necessarily proper) and for any elements $c,d \in F$
if
$$A^{\prime} \subset cG+d,$$
then
$$|A^{\prime}| \leq |G|^{{1 \over 2}}.$$
Then it must be that
$$\max(|A+A|,|AA|) \gtrapprox |A|^{{49 \over 48}}.$$
\endproclaim

The hypotheses regarding interaction with subfields
could be tightened slightly in various ways, but they certainly need to require that $A$ be different
from an affine translate of $G$. In many cases, they are vacuous, as when $k$ is odd and $ p^{{16k+ \epsilon \over 47}}
<|A| < p^{{k \over 2}}$. We think hypotheses as saying that in a certain sense, the dimension of $A$ is at most 
${1 \over 2}$. An analogy may be drawn with the sum product theorem in [B] and it is possible that the techniques
here would be useful in that setting.

\head \S 1 Preliminaries \endhead

Throughout this paper $A$ will denote a fixed set in the field $F=F_{p^k}$ of $p$ elements with $p$
a prime. For $B$, any set, we will denote its cardinality by $|B|$.

Whenever $X$ and $Y$ are quantities we will use
$$X \lesssim Y,$$
to mean 
$$X \leq C Y,$$
where the constant $C$ is universal (i.e. independent of $p$ and $A$). The
constant $C$ may vary from line to line.
We will use
$$X \lessapprox Y,$$
to mean
$$X \leq C (\log |A|)^{\alpha} Y,$$
where $C$ and $\alpha$ may vary from line to line but are universal.

We state some preliminary lemmas.

\proclaim{Lemma 1.1} Let $A \subset F$. Suppose that 
$$|{A-A \over A-A}| \geq |A|^2.$$
Then there are $a_1,a_2,b_1,b_2 \in A$ with
$$|(a_1-a_2) A + (b_1-b_2) A| \gtrsim |A|^2.$$
\endproclaim

\demo{Proof} Under the hypothesis, there is $x \in {A-A \over A-A}$ with at most $|A|^2$
representations
$$x={a_1-a_2 \over b_1-b_2}.$$
Thus there $\lesssim |A|^2$ solutions of
$$a_1 + b_2 x=a_2 + b_1 x.$$
Therefore
$$|A + x A| \gtrsim |A|^2.$$
But if
$$x={a_1-a_2 \over b_1-b_2},$$
Then
$$|A+x A|=|(a_1-a_2) A +(b_1-b_2)A|.$$
\qed \enddemo

\proclaim{Lemma 1.2}  Let $A \subset F$. Suppose that $x \in F$
with $x \notin {A-A \over A-A},$ then
$$|A+xA|=|A|^2.$$
\endproclaim

\demo{Proof}
There are no nontrivial solutions of
$$a_1+x b_2=a_2 + x b_1,$$
with $a_1,a_2,b_1,b_2 \in A$.
\qed \enddemo

\proclaim{Lemma 1.3} Let $A \subset F$ with cardinality at least 3. Suppose that $G$ is a subfield of $F$ with
$${A-A \over A-A} \subset G,$$
then there exist $c,d \in F$ with
$$A \subset cG+d.$$
\endproclaim

\demo{Proof} Suppose that the conclusion is false for all $c,d \in F$. Then we can find $a_1,a_2,b_1,b_2,c,d_1,d_2 \in A$\
and $g_1,g_2,g_3,g_4 \in G$
with $b_1 \neq b_2$ and ${d_1-d_2 \over c} \notin G$, so that
$$a_1= c g_1 + d_1; \quad a_2=c g_2 + d_2; \quad b_1=cg_3 + d_2; \quad b_2=c g_4 + d_2.\tag 1.1$$
We do this as follows: We select $b_1,b_2$ distinct in $A$. Since $b_1-b_2$ is invertible, we can find $c$
so that $(b_1-b_2) \in cG$. Then there is $d_2$ with $b_1,b_2 \in cG+d_2$. We choose $a_2 \in cG+d_2 \cap A$.
It need not be distinct from $b_1$ and $b_2$. Then we apply the assumption to pick $a_1 \in A$ but
$a_1 \notin cG+d_2$. Applying (1.1), we see immediately
$${a_1-a_2 \over b_1-b_2} \notin G.$$
\qed \enddemo

The following two lemmas, quoted by Garaev, are due to Ruzsa, may be found in [TV].
The first is usually referred to as Rusza's triangle inequality. The second is
a form of Plunneke's inequality.

\proclaim{Lemma 1.4} For any subsets $X,Y,Z$ of $F$ we have
$$|X-Z| \leq {|Y-X||X+Z| \over |X|}.$$ \endproclaim

\proclaim{Lemma 1.5} Let $X,B_1,\dots,B_k$ be any subsets of $F$ with
$$|X+B_i| \leq \alpha_i |X|,$$
for $i$ ranging from 1 to $k$. Then there exists $X_1 \subset X$
with
$$|X_1+B_1 + \dots + B_k| \leq \alpha_1 \dots \alpha_k |X_1|. \tag (1.1)$$
\endproclaim

We record a number of Corollaries. The first two can be found in [TV]. The second
one, we first became aware of in the paper of Garaev. The third is a slight refinement
which we need here.

\proclaim{Corollary 1.6} Let $X,B_1,\dots,B_k$ be any subsets of $F$.
Then
$$|B_1 + \dots + B_k| \leq {|X+B_1| \dots |X+B_k| \over |X|^{k-1}}.$$
\endproclaim

\demo{Proof} Simply bound $|B_1 + \dots + B_k|$ by $|X_1+B_1 + \dots + B_k|$
and $|X_1|$ by $|X|$. \qed \enddemo

\proclaim{Corollary 1.7} Let $A \subset F$ and let $a,b \in A$. Then we have
the inequalities
$$|aA+bA| \leq {|A+A|^2 \over |aA \cap bA|},$$
and
$$|aA-bA| \leq {|A+A|^2 \over |aA \cap bA|}.$$
\endproclaim

\demo{Proof}
To get the first inequality, apply Corollary 1.6 with $k=2$, $B_1=aA$, $B_2=bA$, and
$X=aA \cap bA$.

To get the second inequality, apply Lemma 1.4 with $Y=aA$, $Z=-bA$ and $X=-(aA \cap bA)$.
\qed \enddemo

\proclaim{Corollary 1.8} Let $a_1,a_2,b \in F$ and let $A \subset F$. Then
$$|a_1 a_2 A + b^2 A| \leq { |A+A|^4 \over |a_1 A \cap b A| |a_2 A \cap  b A| |A|},$$
and
$$|a_1 a_2 A - b^2 A| \leq { |A+A|^4 \over |a_1 A \cap b A| |a_2 A \cap  b A| |A|},$$
\endproclaim

\demo{Proof}
To obtain the first inequality, we use Corollary 1.6 with $k=2$ and $X= a_1 b A$ to obtain
$$ |a_1 a_2 A + b^2 A| \leq {|a_1 A + b A| |a_2 A + b A| \over |A|}.$$
Then we apply Corollary 1.7 twice. The second inequality proceeds likewise.
\qed \enddemo

\head \S 2 Modified Garaev's inequality  \endhead

In this section, we slightly modify Garaev's argument to obtain the desired result.

\demo{Proof of main Theorem}

Following Garaev, we observe that
$$\sum_{a \in A} \sum_{b \in A} |aA \cap bA| \geq {|A|^4 \over |AA|}.$$
Therefore, we can find an element $b_0 \in A$, a subset $A_1 \subset A$
and a number $N$ satisfying 
$$|b_0 A \cap a A| \approx N,$$
for every $a \in A_1$. Further
$$N \gtrapprox {|A|^2 \over |AA|}, \tag 2.1$$
and 
$$|A_1| N \gtrapprox {|A|^3 \over |AA|}. \tag 2.2 $$

Now there are three cases. In the first case, 
we have that ${A_1 - A_1 \over A_1 - A_1}$ is a field $G \subset F$.
If we have $|A_1| \leq |A|^{{47 \over 48}}$, then we already have the desired
result from (2.2) and $N \leq |A|$. Otherwise, by Lemma 1.3, we have that $A_1$ is
contained in an affine image of $G$ so that by hypothesis
$$|{A_1 - A_1 \over A_1 - A_1}| \gtrsim |A_1|^2.$$
Thus by Lemma 1.1
we can find $a_1,a_2,b_1,b_2 \in A_1$ so that
$$|A_1|^2 \lesssim |(a_1-a_2) A_1 + (b_1-b_2) A_1| \leq |a_1 A - a_2 A + b_1 A - b_2 A|.$$

Applying Corollary 1.6 with $k=4$ and with $B_1= a_1 A$, with $B_2=-a_2A$ with
$B_3 = b_1 A$, with $B_4=-b_2 A$, and with $X=b_0 A$. and applying
Corollary 1.7 to bound above $|X+B_j|$, This we get

$$|A_1|^2 \lesssim {|A+A|^8 \over N^4 |A|^3},$$
or
$$|A_1|^2 N^4 |A|^3 \lessapprox |A+A|^8 .$$
Applying (2.2), we get
$$N^2 |A|^9 \lessapprox |A+A|^8 |AA|^2. \tag 2.3$$
and applying (2.1), we get
$$|A|^{13} \lessapprox |A+A|^8 |AA|^4. \tag 2.4$$
The estimate (2.4) implies that 
$$\max(|A+A|,|AA|) \gtrapprox |A|^{{13 \over 12}}
\gtrapprox |A|^{{49 \over 48}},$$
so that we have more than we need in this case. We
restrict to the setting where ${A_1-A_1 \over A_1 - A_1}$ is not a field.

Now there are two remaining cases, either we can find $x \in {A_1 - A_1 \over A_1 - A_1} +
{A_1 - A_1 \over A_1 - A_1}$ or $x \in ({A_1 - A_1 \over A_1 - A_1})({A_1 - A_1 \over A_1 - A_1})$
with $x \notin {A_1 - A_1 \over A_1 - A_1}$. In light of Lemma 1.2, we then have
$$|A_1|^2 \leq |A_1 + x A_1| \leq |A+xA|.$$

In the former case, we may write $x={a_1 - a_2 \over b_1-b_2} + {c_1-c_2 \over d_1-d_2}$.
Then we have
$$|A_1|^2 \leq |(b_1-b_2)(d_1-d_2) A + (b_1-b_2)(c_1-c_2) A + (a_1-a_2)(d_1-d_2) A|.$$
Now applying Corollary 1.6 with $k=3$  and $X=(b_1-b_2)(d_1-d_2) A$, we obtain
$$|A_1|^2 \leq {|A+A| |d_1 A-d_2  A + c_1 A-c_2 A| |b_1 A-b_2 A + d_1 A-d_2 A| \over |A|^2}.$$
Applying Corollary 1.6 to the last two factors with $k=4$ and $X=b_0 A$ and then invoking
Corollary 1.7, we get
$$|A_1|^2 \lesssim {|A+A|^{17} \over N^8 |A|^8}.$$
Invoking (2.1) and (2.2), we get
$$|A|^{26} \lessapprox |A+A|^{17} |AA|^8.$$
Since ${26 \over 25} \geq {49 \over 48},$ we get more than we need in this case as well.

In the final case, we may write
$$x={(a_{10}-a_{11})(a_{20}-a_{21}) \over (a_{30}-a_{31})(a_{40}-a_{41})}.$$
Thus we obtain
$$|A_1|^2 \leq | \sum_{i=0}^1 \sum_{j=0}^1  (-1)^{i+j} a_{1i} a_{2j} A +
\sum_{i=0}^1 \sum_{j=0}^1  (-1)^{i+j} a_{3i} a_{4j} A|.$$
We now apply Corollary 1.6 with $k=8$ and $X=b_0^2 A$ and we apply Corollary 1.8 to each of the
eight terms to obtain
$$|A_1|^2 \leq {|A+A|^{32} \over N^{16} |A|^{15} }.$$
Applying (2.1) and (2.2), we obtain
$$|A|^{49} \lessapprox  |A+A|^{32} |AA|^{16},$$
which gives the desired result.

\qed \enddemo

\end